\documentclass[11pt]{article}

\usepackage{amsmath}
\usepackage{amssymb}
\usepackage{authblk}
\usepackage{booktabs}
\usepackage{color}
\usepackage{epsfig,subfigure,stmaryrd,appendix}
\usepackage{epstopdf}
\usepackage{graphicx}
\usepackage{latexsym}
\usepackage{amsfonts}

\begin{document}  \baselineskip=16pt

\centerline{\LARGE Superconvergence Points of Spectral Interpolation}

\bigskip
\centerline{Zhimin Zhang \footnote{Department of
Mathematics, Wayne State University, Detroit, MI 48202, USA.
This work is supported in part by
the US National Science Foundation under grant DMS-1115530.}}

\bigskip


{\it ABSTRACT. In this work, we study superconvergence properties for
some high-order orthogonal polynomial interpolations.
The results are two-folds: When interpolating function values,
we identify those points where the first and second derivatives of the interpolant converge faster;
When interpolating the first derivative,
we locate those points where the function value of the interpolant superconverges.
For the earlier case, we use various Chebyshev polynomials; and for the later case,
we also include the counterpart Legendre polynomials.}




\bigskip
\setcounter{equation}{0}
\setcounter{section}{1}

\noindent{\bf 1. Introduction}

In numerical computation, we often observe that the convergent rate exceeds
the best possible global rate at some special points.
Those points are called superconvergent points,
and the phenomenon is called superconvergence phenomenon, which
is well understood for the $h$-version finite element method, see, e.g.,
\cite{babuska-sug, chen, chen-huang, krizek-n, krizek-ns, lin-lin, lin-yan, schatz-sloan-wahlbin, wahlbin, zhu-lin}
and references therein.
As comparison, the relevant study for the $p$-version finite element method and the spectral method is lacking.
Only very special and simple cases have been discussed in the following works:
 \cite{wang-xie-zhang, zhang1, zhang2}.

The study of the superconvergence phenomenon for the $h$-version method has made great impact
on scientific computing, especially on {\it a posteriori} error estimates and adaptive methods,
which is well documented in the following books:
\cite{ainsworth-oden, babuska-strouboulis, krizek-ns, tang-xu} and works cited.
It is the believe of this author that the scientific community would also be benefited from the
study of superconvergence phenomenon of spectral collocation methods as well as related $p$-version
and spectral methods. This work is the first step, where the superconvergence points
of some orthogonal polynomial interpolation will be identified.

The most celebrated advantage of spectral methods is the exponential
(or geometric) rate of convergence
for sufficiently smooth, essentially analytic functions.
However, most error bounds in the literature
are in the form of $N^{-k}\|u\|_{H^{k+1}}$, where $N$ is the polynomial degree
(or trigonometric function degree in the Fourier spectral case).
We can see this in almost all books on spectral methods such as
\cite{bernardi-maday, boyd, canuto-hqz, canuto-hqz1, canuto-hqz2, deville-fm, fornberg, funaro,
gottlieb-orszag, guo, hesthaven-gg, karniadakis-sherwin, peyret, shen-tang, shen-tang-wang}
Ideally, we expect to establish the convergence rate $\rho^{-N}$ for some $\rho>1$
or $e^{-\sigma N}$ for some $\sigma>0$.
There have been some limited discussion of this type of error bounds in the past,
 e.g., in these two books \cite{gottlieb-orszag, trefethen},
 see also \cite{gottlieb-shu, reddy-weideman, tadmor}.
In the framework of $p$- and $hp$-finite element methods,
results of the exponential convergent rate
can be found in books \cite{schwab, szabo-babuska} and references therein,
see also, \cite{guo-babuska1, guo-babuska2}.

Actually, exponential rate of convergence for polynomial approximation of analytic functions
can be traced back 100 years. The following result is due to Bernstein \cite{bernstein},
see also \cite[Theorem 7]{lorentz}: $f$ is analytic on $[-1,1]$ if and only if
$$
\sup\lim_{n\to\infty} \sqrt[n]{E_n(f)} = \rho^{-1}, \quad
E_n(f) = \inf_{g\in{\cal P}_n} \|f-g\|_\infty,
$$
where ${\cal P}_n$ is the polynomial space of degree no more than $n$,
and $\rho>1$ is the sum of the half-axes of the maximum elliptic disc $D_\rho$
bounded by the ellipse $E_\rho$ with foci $\pm 1$ that $f$ can be analytically extended to.
This fact serves as a starting point of the first part of our analysis in this paper.
We shall focus our attention to the approximation properties of
polynomial interpolants of analytic functions on $[-1,1]$.
We identify those points in $[-1,1]$, where
the derivatives are superconvergent, in the sense that the convergent rate gains
at least one factor $N^{-1}$. Three popular Chebyshev polynomial interpolants
are discussed here: Chebyshev, Chebyshev-Lobatto, and Chebyshev-Radau.

In the second part of analysis, we consider polynomial interpolants of the first derivative
of a one-order higher polynomial, and identify those points where function values
are superconvergent. In this part, both the Chebyshev and Legendre polynomials
are studied.

Throughout of the paper, we use the standard notation for orthogonal polynomials:
$T_j$, the Chebyshev polynomial of the first kind;
$U_j$, the Chebyshev polynomial of the second kind; and
$L_j$, the Legendre polynomial.

\bigskip


\setcounter{equation}{0}
\setcounter{section}{2}

\noindent{\bf 2. Interpolation by the Chebyshev Polynomials}

\vskip.05in

In this Section, we discuss interpolations by the Chebyshev polynomials of the first kind,
Chebyshev-Lobatto polynomials, and Chebyshev-Radau (right and left) polynomials.

\vskip.05in

\noindent {\bf 2.1. Statement of the Results}

\vskip.05in

We state the superconvergence results in this subsection.
The meaning of the superconvergence will be made clear in the next subsection.

\noindent {\bf 2.1.1. Chebyshev interpolant}. Interpolating at the zeros of $T_{N+1}$,
\begin{equation}\label{cheb0}
T_{N+1}(x_k) = \cos(N+1)\theta_k = 0, \quad x_k = \cos\theta_k, \quad
\theta_k = \frac{2k+1}{2N+2}\pi, \quad k=0,1,\ldots,N.
\end{equation}
{\bf Proposition 2.1}. For the Chebyshev interpolant, the first derivative superconverges
at zeros of $U_N$, which are
\begin{equation}\label{cheb1}
y_k = \cos\frac{k\pi}{{N+1}}, \quad k=1,2,\ldots,N;
\end{equation}
and the second derivative superconverges at $\cos\theta_k$ with $\theta_k$ satisfies the
following equation:
\begin{equation}\label{cheb2}
(N+1)\cos(N+1)\theta\sin\theta = \sin(N+1)\theta\cos\theta.
\end{equation}
This set of $\cos\theta_k$ are close to interior zeros
(not including $k=0,N$ in (\ref{cheb0})) of $T_{N+1}$ for large $N$.

\vskip.1in

\noindent {\bf 2.1.2. Chebyshev-Lobatto interpolant}.
Interpolating at the zeros of $T_{N+1}-T_{N-1}$,
\begin{equation}\label{cheb-loba}
(T_{N+1}-T_{N-1})(x_k) = 2\sin N\theta_k\sin\theta_k = 0, \quad
x_k = \cos\frac{k\pi}{N}, \quad k=0,1,\ldots,N.
\end{equation}
We also call this set of points the 2nd type Chebyshev points.

\noindent {\bf Proposition 2.2}. For the Chebyshev-Lobatto interpolant,
the first derivative superconverges at
$\cos\theta_k$ with $\theta_k$ satisfies the following equation:
\begin{equation}\label{cheb-loba1}
N\cos N\theta\sin\theta + \sin N\theta\cos\theta = 0.
\end{equation}
This set of $\theta_k$s are close to zeros of $\cos N\theta$ for large $N$,
i.e., $\displaystyle \theta_k \approx \frac{2k-1}{2N}\pi$ when $\theta_k$ is away from $0$ and $\pi$.

As for superconvergent points of the second derivative,
$\theta_k$'s satisfy the following equation:
\begin{equation}\label{cheb-loba2}
(N^2-1)\sin N\theta \sin\theta = 2N \cos N\theta \cos\theta.
\end{equation}
This set of $\theta_k$s are close to interior zeros (not including 0 and $\pi$)
of $\sin N\theta$ for large $N$,
i.e., $\displaystyle \theta_k \approx \frac{k\pi}{N}$, when $\theta_k$ is away from $0$ and $\pi$.

\vskip.1in

{\it Remark 2.1}. We see that when interpolating at the zeros of $T_{N+1}(x)$,
the derivative of the interpolant superconverges at the zeros of $U_N(x)$,
and the second derivative superconvergent points go back ``almost"
to the zeros of $T_{N+1}(x)$ except the two ends;
when interpolating at the zeros of $(T_{N+1}-T_{N-1})(x) = \gamma_N(1-x^2)U_{N-1}(x)$,
the derivative of the interpolant superconverges ``almost" at the zeros of $T_N(x)$,
and the second derivative superconvergent points go back ``almost"
to the zeros of $U_{N-1}(x)$.

\vskip.1in

\noindent {\bf 2.1.3. Chebyshev-Radau (right) interpolant}.
Interpolating at the zeros of $T_{N+1}-T_N$,
\begin{equation}\label{cheb-right}
(T_{N+1}-T_N)(x_k) = 2\sin(N+\frac{1}{2})\theta_k \sin\frac{\theta_k}{2} = 0, \quad
x_k = \cos\frac{2k\pi}{2N+1}, \quad k=0,1,\ldots,N.
\end{equation}
{\bf Proposition 2.3}. For the right Chebyshev-Radau interpolant,
the first derivative superconverges at
$\cos\theta_k$ with $\theta_k$ satisfies the following equation:
\begin{equation}\label{cheb-right1}
(2N+1)\cos(N+\frac{1}{2})\theta \sin\frac{\theta}{2}
+ \sin(N+\frac{1}{2})\theta \cos\frac{\theta}{2} = 0.
\end{equation}
This set of $\theta_k$s are close to zeros of $\cos(N+\frac{1}{2})\theta$ for large $N$,
i.e., $\displaystyle \theta_k \approx \frac{2k-1}{2N+1}\pi$, when $\theta_k$ is away from $0$.

As for superconvergent points of the second derivative,
$\theta_k$ satisfies the following equation:
\begin{equation}\label{cheb-right2}
(2N^2+2N+1)\sin(N+\frac{1}{2})\theta \sin\frac{\theta}{2}
= (2N+1)\cos(N+\frac{1}{2})\theta \cos\frac{\theta}{2}.
\end{equation}
This set of $\theta_k$s are close to interior zeros (not including 0)
of $\sin(N+\frac{1}{2})\theta$ for large $N$,
i.e., $\displaystyle \theta_k \approx \frac{2k\pi}{2N+1}$, when $\theta_k$ is away from $0$.

\vskip.1in

\noindent {\bf 2.1.4. Chebyshev-Radau (left) interpolant}.
Interpolating at the zeros of $T_{N+1}+T_N$,
\begin{equation}\label{cheb-left}
(T_{N+1}+T_N)(x_k) = 2\cos(N+\frac{1}{2})\theta_k \cos\frac{\theta_k}{2} = 0, \quad
x_k = \cos\frac{2k+1}{2N+1}\pi, \quad k=0,1,\ldots,N.
\end{equation}
{\bf Proposition 2.4}. For the left Chebyshev-Radau interpolant,
the first derivative superconverges at
$\cos\theta_k$ with $\theta_k$ satisfies the following equation:
\begin{equation}\label{cheb-left1}
(2N+1)\sin(N+\frac{1}{2})\theta \cos\frac{\theta}{2}
+ \cos(N+\frac{1}{2})\theta \sin\frac{\theta}{2} = 0.
\end{equation}
This set of $\theta_k$s are close to zeros of $\displaystyle\sin(N+\frac{1}{2})\theta$ for large $N$,
i.e., $\displaystyle \theta_k \approx \frac{2k}{2N+1}\pi$, when $\theta_k$ is away from $\pi$.

As for superconvergent points of the second derivative,
$\theta_k$ satisfies the following equation:
\begin{equation}\label{cheb-left2}
(2N^2+2N+1)\cos(N+\frac{1}{2})\theta \cos\frac{\theta}{2}
= (2N+1)\sin(N+\frac{1}{2})\theta \sin\frac{\theta}{2}.
\end{equation}
This set of $\theta_k$s are close to interior zeros (not including $\pi$)
of $\cos(N+\frac{1}{2})\theta$ for large $N$,
i.e., $\displaystyle \theta_k \approx \frac{2k+1}{2N+1}\pi$, when $\theta_k$ is away from $\pi$.

\vskip.1in

{\it Remark 2.2}. Propositions 2.3 and 2.4 say that when interpolating at the right (left) Radau points,
derivative of the interpolant superconverges ``almost" at the left (right) Radau points except $x=-1$ ($x=1$),
and the second derivative superconvergent points ``almost" go back
to the interior right (left) Radau points.

{\it Remark 2.3}. In this section, we provide derivative superconvergence points as roots of some polynomial equations, and the approximated values of those points. Numerical data in the last section demonstrate that these approximation are usually good enough for practical purpose. Therefore, they can be used as the intitial guesses, e.g., the Newton iteration, should the more accurate values are desired. On the other hand, there are quick and simple ways to find roots of sums of orthogonal polynomials by expressing the problem as an eigenvalue problem wuith the matrix being upper Hessenberg, see \cite{day-romero}.

\vskip.1in

\noindent{\bf 2.2. Analysis}

\vskip.05in

Let $u$ be analytic on $I=[-1,1]$. According to Bernstein \cite{bernstein},
$u$ can be analytically extended to $D_\rho$, which is
enclosed by an ellipse $E_\rho$ with $\pm 1$ as fosi,
$\rho>1$ as the sum of its semimajor and semiminor:
$$
E_\rho : \qquad z = \frac{1}{2} ( \rho e^{i\theta} + \rho^{-1}e^{-i\theta} ),
\qquad  0\le \theta \le 2\pi.
$$
We consider polynomial $u_N \in {\cal P}_N$ who interpolates $u$ at $N+1$ points
$-1 \le x_0 \le x_1 \le \cdots \le x_N \le 1$.
The error equation is, according to \cite[p.68]{davis}, expressed as
\begin{equation}\label{err-eq0}
u(x) - u_N(x) = \frac{1}{2\pi i} \int_{E_\rho} \frac{\omega_{N+1}(x)}{\omega_{N+1}(z)}
\frac{u(z)}{z-x} dz, \quad
\omega_{N+1}(x) = c\Pi_{j=0}^N (x-x_j).
\end{equation}
By direct differentiation of (\ref{err-eq0}), c.f., \cite{reddy-weideman},
one obtains the error equation for the derivative
\begin{equation}\label{err-eq1}
u'(x) - u_N'(x) = \frac{1}{2\pi i} \int_{E_\rho} \left(
\frac{\omega_{N+1}'(x)}{z-x} + \frac{\omega_{N+1}(x)}{(z-x)^2} \right)
\frac{u(z)}{\omega_{N+1}(z)} dz,
\end{equation}
and the second derivative
\begin{equation}\label{err-eq2}
u''(x) - u_N''(x) = \frac{1}{2\pi i} \int_{E_\rho} \left(
\frac{\omega_{N+1}''(x)}{z-x} + \frac{2\omega_{N+1}'(x)}{(z-x)^2}
+ \frac{2\omega_{N+1}(x)}{(z-x)^3} \right)
\frac{u(z)}{\omega_{N+1}(z)} dz.
\end{equation}

We need to estimate $\omega_{N+1}(x)$, $\omega_{N+1}'(x)$,
and $\omega_{N+1}''(x)$ on $[-1,1]$ to establish the error bounds.
However, by each differentiation, we lose at least one power of $N$.

{\bf Key observation}. Let us examine the error equation (\ref{err-eq1}).
At the $N$ special points $\omega_{N+1}'(x) = 0$, we have only the
second term, which is usually smaller than the first term in magnitude by
a factor $N$ or $N^2$ as we will see later. Similarly, at the
$N-1$ special points $\omega_{N+1}''(x) = 0$, we have only the
second and third terms left in the error equation (\ref{err-eq2}).
Again, we may gain a factor $N$ in the error bounds.

We consider the four sets of interpolation points in the previous subsection.

The exponential decay of the error is provided by the value of $\omega_{N+1}(z)$
on the ellipse $E_\rho$ in the denominators of (\ref{err-eq0})-(\ref{err-eq2}).
In all four sets of interpolation points, $\omega_{N+1}(z)$ involves the Chebyshev
polynomials of the first kind. We have the following characteristic expressions.
The proof is elementary by the definition of $T_{N+1}(z)$ and therefore is omitted.

\vskip.1in

\noindent {\bf Lemma 2.1}. For $z\in E_\rho$, we have
\begin{equation}
\label{ellp}
| T_{N+1}(z) | = \frac{1}{2} \sqrt{ \rho^{2N+2} + \rho^{-2N-2} + 2\cos 2(N+1)\theta };
\end{equation}
\begin{equation}
\label{ellp-loba}
| T_{N+1}(z) - T_{N-1}(z) | = \frac{1}{2} \sqrt{ \rho^2 + \rho^{-2} - 2\cos2\theta }
\sqrt{ \rho^{2N} + \rho^{-2N} - 2\cos2N\theta };
\end{equation}
\begin{equation}
\label{ellp-radau}
| T_{N+1}(z) \pm T_N(z) | = \frac{1}{2} \sqrt{ 1 + \rho^{-2} \pm 2\rho^{-1}\cos\theta }
\sqrt{ \rho^{2N+1} + \rho^{-2N} \pm 2\rho\cos(2N+1)\theta }.
\end{equation}

We also need to bound the derivatives of $\omega_{N+1}(x)$ for $x\in [-1,1]$.

\vskip.1in

\noindent {\bf Lemma 2.2}.
\begin{eqnarray}
\label{d1}
\max_{x\in[-1,1]} | T_{N+1}'(x) | &=& (N+1)^2, \\
\label{d1-loba}
\max_{x\in[-1,1]} | T_{N+1}'(x) - T_{N-1}'(x) | &=& 4N, \\
\label{d1-radau}
\max_{x\in[-1,1]} | T_{N+1}'(x) \pm T_N'(x) | &=& 2N^2+2N+1; \\
\label{d2}
\max_{x\in[-1,1]} | T_{N+1}''(x) | &=& \frac{1}{3}N(N+1)^2(N+2), \\
\label{d2-loba}
\max_{x\in[-1,1]} |T_{N+1}''(x)-T_{N-1}''(x) | &=& \frac{4}{3} N(2N^2+1), \\
\label{d2-radau}
\max_{x\in[-1,1]} | T_{N+1}''(x) \pm T_N''(x) | &=& \frac{2}{3} N(N+1)(N^2+N+1).
\end{eqnarray}

The proof is also elementary and is provided in the appendix.
To establish the error bounds, we define some constants,
which include $C_\rho(u) = \displaystyle \max_{z\in E_\rho} |u(z)|$,
$D_\rho$, the shortest distance from $E_\rho$ to $[-1,1]$,
and $L_\rho$, the arch length of the ellipse $E_\rho$. We have
\begin{equation}\label{err7}
D_\rho = \frac{1}{2} (\rho+\rho^{-1}) - 1, \qquad L_\rho \le \pi\sqrt{\rho^2+\rho^{-2}}.
\end{equation}
The latter is the Euler's estimate which overestimates the perimeter by less than 12 percent.

\noindent {\bf Theorem 2.1}. Assume that $u$ is analytic on $[-1,1]$ and can be extended analytically
to the region bounded by an ellipse $E_\rho$. Let $u_N\in {\cal P}_N[-1,1]$ be the interpolant of $u$
at $N+1$ zeros of $T_{N+1}(x)$ on $[-1,1]$. Then
\begin{equation}\label{err-t0}
\max_{x\in[-1,1]} |u(x) - u_N(x) | \le \frac{C_\rho(u) L_\rho}{\pi D_\rho} \frac{1}{\rho^{N+1}-\rho^{-N-1}},
\end{equation}
\begin{equation}\label{err-t1}
\max_{x\in[-1,1]} |u'(x)-u_N'(x) | \le \frac{C_\rho(u) L_\rho}{\pi D_\rho}
\left( (N+1)^2 + \frac{1}{D_\rho} \right) \frac{1}{\rho^{N+1}-\rho^{-N-1}},
\end{equation}
\begin{equation}\label{err-t2}
\max_{x\in[-1,1]} |u''(x) - u_N''(x)| \le \frac{C_\rho(u) L_\rho}{\pi}
\left( \frac{N(N+1)^2(N+2)}{3D_\rho} + \frac{2(N+1)^2}{D_\rho^2} + \frac{2}{D_\rho^3}\right) \frac{1}{\rho^{N+1}-\rho^{-N-1}};
\end{equation}
Furthermore, at those special points where $T_{N+1}'(x) = 0$, we have
\begin{equation}\label{derr-t1}
\max_{1\le j\le N} |u'(t_j)-u_N'(t_j) | \le \frac{C_\rho(u) L_\rho}{\pi D_\rho^2}
\frac{1}{\rho^{N+1}-\rho^{-N-1}}, \quad t_j = \cos\frac{j\pi}{N+1}, \quad j=1,2,\ldots,N,
\end{equation}
and at those special points where $T_{N+1}''(x) = 0$, we have
\begin{equation}\label{derr-t2}
\max_{1\le j\le N-1} |u''(\tau_j) - u_N''(\tau_j)| \le \frac {2C_\rho(u) L_\rho}{\pi D_\rho^2}
\left( (N+1)^2 + \frac{1}{D_\rho} \right) \frac{1}{\rho^{N+1}-\rho^{-N-1}},
\end{equation}
where $\tau_j = \cos\theta_j$ with $\theta_j$ satisfies (\ref{cheb2}).

Proof: Since $u_N\in {\cal P}_N[-1,1]$ interpolates $u$
at zeros of $T_{N+1}(x)$, we have $\omega_{N+1} = T_{N+1}$ in
(\ref{err-eq0}) - (\ref{err-eq2}).
By the identity (\ref{ellp}), we can derive the lower bound
$$
|T_{N+1}(z)| \ge \frac{1}{2} (\rho^{N+1}-\rho^{-N-1}).
$$
Substituting this into (\ref{err-eq0}) and using
$\displaystyle \max_{x\in[-1,1]} |T_{N+1}(x)| = 1$, we derive,
\begin{equation}\label{t3.1}
| u(x) - u_N(x)| \le \frac{1}{\pi} \int_{E_\rho} \frac{|u(z)|}{|z-x|} d|z|
\frac{1}{\rho^{N+1}-\rho^{-N-1}}, \quad \forall x\in[-1,1].
\end{equation}
Note that $|u(z)| \le C_\rho(u)$
and $|z-x|^{-1} \le D_\rho^{-1}$ for $z\in E_\rho$ and $x\in[-1,1]$, and
we obtain (\ref{err-t0}) from (\ref{t3.1}).

Using (\ref{d1}) in (\ref{err-eq1}) and following the same procedure
as above, we derive (\ref{err-t1}). Similarly, we establish (\ref{err-t2})
by applying (\ref{d1}) and (\ref{d2}) in (\ref{err-eq2}).

At the special points when $T_{N+1}'(x) = (N+1)U_N(x) = 0$, the first term on the right-hand side
of (\ref{err-eq1}) is gone, we then obtain (\ref{derr-t1}) following the same
argument as we derive (\ref{err-t0}).

At the special points when $T_{N+1}''(x) = (N+1)U_N'(x) = 0$, the first term on the right-hand side
of (\ref{err-eq2}) is gone, we then obtain (\ref{derr-t2}) following the same
argument as we derive (\ref{err-t1}). Since $U_N(x)$ has $N$ simple roots in $(-1,1)$,
it is guaranteed that $U_N'(x)$ has $N-1$ simple roots in between and they can be
expressed as $x_k = \cos\theta_k$ with $\theta_k$ satisfying (\ref{cheb2}).

\vskip.2in

\noindent {\bf Theorem 2.2}. Under the same assumption as in Theorem 2.1,
let $u_N\in {\cal P}_N[-1,1]$ be the interpolant of $u$
at $N+1$ zeros of $T_{N+1}(x)-T_{N-1}(x)$ on $[-1,1]$. Then
\begin{equation}\label{err-l0}
\max_{x\in[-1,1]} |u(x) - u_N(x) | \le \frac{C_\rho(u) L_\rho}{\pi D_\rho} \frac{(\rho-\rho^{-1})^{-1}}{\rho^N-\rho^{-N}},
\end{equation}
\begin{equation}\label{err-l1}
\max_{x\in[-1,1]} |u'(x)-u_N'(x) | \le \frac{C_\rho(u) L_\rho}{\pi D_\rho}
\left( 4N + \frac{1}{D_\rho} \right) \frac{(\rho-\rho^{-1})^{-1}}{\rho^N-\rho^{-N}},
\end{equation}
\begin{equation}\label{err-l2}
\max_{x\in[-1,1]} |u''(x) - u_N''(x)| \le \frac{C_\rho(u) L_\rho}{\pi}
\left( \frac{4N(2N^2+1)}{3D_\rho} + \frac{8N}{D_\rho^2} + \frac{2}{D_\rho^3}\right) \frac{(\rho-\rho^{-1})^{-1}}{\rho^N-\rho^{-N}};
\end{equation}
Furthermore, at those special points where $T_{N+1}'(x)-T_{N-1}'(x) = 0$, we have
\begin{equation}\label{derr-l1}
\max_{1\le j\le N} |u'(t_j)-u_N'(t_j) | \le \frac{C_\rho(u) L_\rho}{\pi D_\rho^2}
\frac{(\rho-\rho^{-1})^{-1}}{\rho^N-\rho^{-N}},
\end{equation}
where $t_j = \cos\theta_j$ with $\theta_j$ satisfies (\ref{cheb-loba1});
and at those special points where $T_{N+1}''(x)-T_{N-1}''(x)= 0$, we have
\begin{equation}\label{derr-l2}
\max_{1\le j\le N-1} |u''(\tau_j) - u_N''(\tau_j)| \le \frac {2C_\rho(u) L_\rho}{\pi D_\rho^2}
\left( 4N + \frac{1}{D_\rho} \right) \frac{(\rho-\rho^{-1})^{-1}}{\rho^N-\rho^{-N}},
\end{equation}
where $\tau_j = \cos\theta_j$ with $\theta_j$ satisfies (\ref{cheb-loba2}).

Proof: Since $u_N\in {\cal P}_N[-1,1]$ interpolates $u$
at zeros of $(T_{N+1}-T_{N-1})(x)$, we have $\omega_{N+1} = T_{N+1}-T_{N-1}$ in
(\ref{err-eq0}) - (\ref{err-eq2}).
By the identity (\ref{ellp-loba}), we can derive the lower bound
$$
|(T_{N+1}-T_{N-1})(z)| \ge \frac{1}{2} (\rho-\rho^{-1})(\rho^N-\rho^{-N}).
$$
Using (\ref{d1-loba}) and (\ref{d2-loba}), the rest is then similar to the proof of Theorem 2.1.

\vskip.2in

\noindent {\bf Theorem 2.3}. Under the same assumption as in Theorem 2.1,
let $u_N\in {\cal P}_N[-1,1]$ be the interpolant of $u$
at $N+1$ zeros of $T_{N+1}(x)\pm T_N(x)$ on $[-1,1]$. Then
\begin{equation}\label{err-r0}
\max_{x\in[-1,1]} |u(x) - u_N(x) | \le \frac{C_\rho(u) L_\rho}{\pi D_\rho} \frac{(1-\rho^{-1})^{-1}}{\rho^N-\rho^{-N}},
\end{equation}
\begin{equation}\label{err-r1}
\max_{x\in[-1,1]} |u'(x)-u_N'(x) | \le \frac{C_\rho(u) L_\rho}{\pi D_\rho}
\left( 4N + \frac{1}{D_\rho} \right) \frac{(1-\rho^{-1})^{-1}}{\rho^N-\rho^{-N}},
\end{equation}
\begin{equation}\label{err-r2}
\max_{x\in[-1,1]} |u''(x) - u_N''(x)| \le \frac{C_\rho(u) L_\rho}{\pi}
\left( \frac{4N(2N^2+1)}{3D_\rho} + \frac{8N}{D_\rho^2} + \frac{2}{D_\rho^3}\right) \frac{(1-\rho^{-1})^{-1}}{\rho^N-\rho^{-N}};
\end{equation}
Furthermore, at those special points where $T_{N+1}'(x)\pm T_N'(x) = 0$, we have
\begin{equation}\label{derr-r1}
\max_{1\le j\le N} |u'(t_j)-u_N'(t_j) | \le \frac{C_\rho(u) L_\rho}{\pi D_\rho^2}
\frac{(1-\rho^{-1})^{-1}}{\rho^N-\rho^{-N}},
\end{equation}
where $t_j = \cos\theta_j$ with $\theta_j$ satisfies (\ref{cheb-left1}) in case $T_{N+1}+T_N$, and
(\ref{cheb-right1}) for $T_{N+1}-T_N$;
and at those special points where $T_{N+1}''(x)\pm T_N''(x)= 0$, we have
\begin{equation}\label{derr-r2}
\max_{1\le j\le N-1} |u''(\tau_j) - u_N''(\tau_j)| \le \frac {2C_\rho(u) L_\rho}{\pi D_\rho^2}
\left( 4N + \frac{1}{D_\rho} \right) \frac{(1-\rho^{-1})^{-1}}{\rho^N-\rho^{-N}},
\end{equation}
where $\tau_j = \cos\theta_j$ with $\theta_j$ satisfies (\ref{cheb-left2})
in case $T_{N+1}+T_N$, and (\ref{cheb-right2}) for $T_{N+1}-T_N$.

Proof: Since $u_N\in {\cal P}_N[-1,1]$ interpolates $u$
at zeros of $(T_{N+1}\pm T_N)(x)$, we have $\omega_{N+1} = T_{N+1}\pm T_N$ in
(\ref{err-eq0}) - (\ref{err-eq2}).
By the identity (\ref{ellp-radau}), we can derive the lower bound
$$
|(T_{N+1}\pm T_N)(z)| \ge \frac{1}{2} (1-\rho^{-1})(\rho^N-\rho^{-N}).
$$
Using (\ref{d1-radau}) and (\ref{d2-radau}), the rest is then similar to the proof of Theorem 2.1.

\vskip.1in

We see from the above that superconvergent points involve extremals of $\omega_{N+1}(x)$ for
$x\in[-1,1]$. It would be interesting to see the distribution and magnitudes of those points
for different cases. Since those information is clear for the case $\omega_{N+1}(x) = T_{N+1}(x)$,
we consider the other three cases.

\vskip.1in

\begin{figure} 
        \centering
        \includegraphics[scale=0.5]{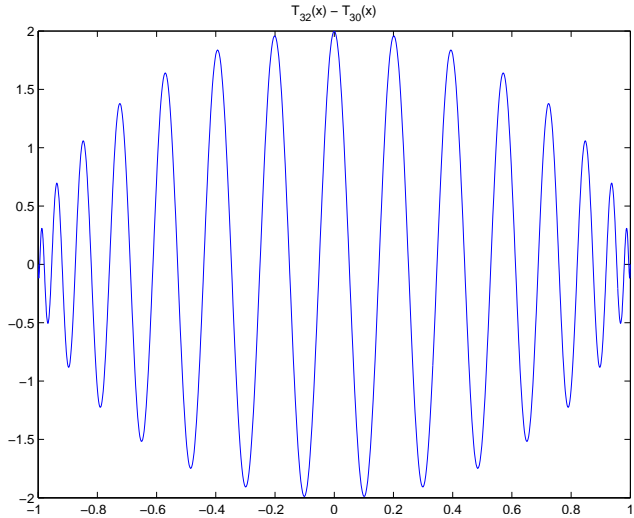}
    \caption{Profile of the Chebshev-Lobatto Polynomial}
    \label{ch-lobatto}
\end{figure}

\noindent {\bf Theorem 2.4}. The envelope for the extremals of $T_{N+1}-T_{N-1}$ on $[-1,1]$
forms an ellipse $\displaystyle x^2 + \frac{y^2}{2^2} = 1$ (see Figure \ref{ch-lobatto}).

Proof: We need to demonstrate that the envelope of $(x,T_{N+1}(x)-T_{N-1}(x))$, or
$$
(\cos\theta, \cos(N+1)\theta - \cos(N-1)\theta) \quad\text{or}\quad
(\cos\theta, -2\sin N\theta\sin\theta)
$$
is an ellipse.
Note that $\sin N\theta$ has extremals at $\displaystyle\theta_j = \frac{2j-1}{2N}\pi$
and points $(\cos\theta_j, \pm 2\sin\theta_j)$ are on the indicated ellipse.

\vskip.1in

\begin{figure} 
        \centering
        \includegraphics[scale=0.5]{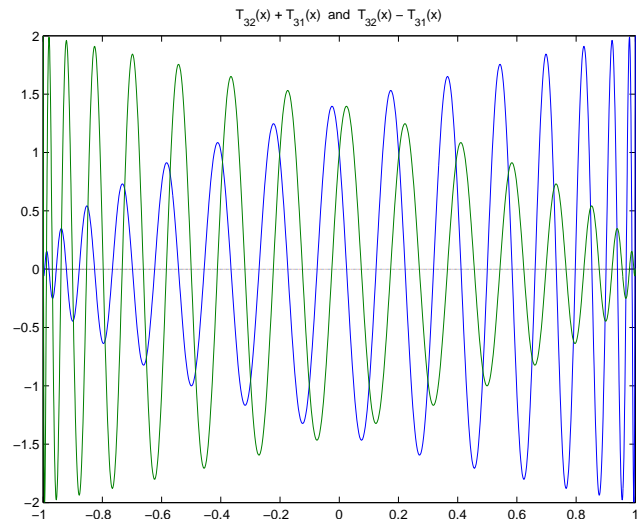}
    \caption{Profile of the Chebshev-Radau Polynomials}
    \label{ch-radau}
\end{figure}

\noindent {\bf Theorem 2.5}. The envelope for the extremals of $T_{N+1} \pm T_N$ on $[-1,1]$
form a parabola $2(1\pm x) = y^2$ (see Figure \ref{ch-radau}).

Proof: We need to demonstrate that the envelope of $(x,T_{N+1}(x)\pm T_N(x))$, or
$$(\cos\theta, \cos(N+1)\theta \pm \cos N\theta ),$$
or
$$
(\cos\theta, 2\cos(N+\frac{1}{2})\theta\cos\frac{\theta}{2}) \quad \text{and}\quad
(\cos\theta, - 2\sin(N+\frac{1}{2})\theta\sin\frac{\theta}{2})
$$
are parabola. Note that
$$
\cos(N+\frac{1}{2})\theta \quad\text{and}\quad \sin(N+\frac{1}{2})\theta
$$
have extremals at
$$
\theta_j = \frac{2j\pi}{2N+1} \quad\text{and}\quad \theta_j = \frac{2j+1}{2N+1}\pi,
$$
respectively, and points
$$
(\cos\theta_j, \pm 2\cos\frac{\theta_j}{2}) \quad\text{and}\quad
(\cos\theta_j, \pm 2\sin\frac{\theta_j}{2})
$$
are on the indicated parabola, respectively. Indeed,
$$
2(1+\cos\theta) = (2\cos\frac{\theta}{2})^2  \quad\text{and}\quad
2(1-\cos\theta) = (2\sin\frac{\theta}{2})^2.
$$

\bigskip


\setcounter{equation}{0}
\setcounter{section}{3}

\noindent{\bf 3. Derivative Interpolation}

In this part, we consider different interpolants to the first derivative of a smooth function,
and identify superconvergent points for the function value approximation.
To be more precise, we construct polynomial $u_N \in{\cal P}_N$ such that
\begin{equation}\label{ode}
u_N(-1) = u(-1), \quad u_N'(x_k) = u'(x_k), \quad k=1,2,\ldots,N,
\quad -1 \le x_1 < \cdots < x_N \le 1,
\end{equation}
and locate the point $y_k$, where $(u-u_N)(y_k)$ is superconvergent.
It is worthy to point out that this supeconvergence knowledge can be utilized in
spectral collocation method for solving ODEs. We shall demonstrate this point later
in our numerical examples.

To fix the idea, we consider only the case $u\in {\cal P}_{N+1}$, since superconvergent property
may be narrowed down to the capability of a polynomial space to approximate polynomials of one order higher
\cite{babuska-sug}.

In addition to Chebyshev, Chebyshev-Lobatto, and Chebyshev-Radau interpolants,
we also consider Gauss, Gauss-Lobatto, and Gauss-Radau interpolants.

We begin with two technical lemmas.

\vskip.1in

\noindent{\bf Lemma 3.1}.
\begin{equation}\label{radau1}
N(L_N+L_{N-1})(x) = (x+1)(L_N-L_{N-1})'(x).
\end{equation}
\begin{equation}\label{radau2}
N(L_N-L_{N-1})(x) = (x-1)(L_N+L_{N-1})'(x).
\end{equation}

Proof: We only prove (\ref{radau1}). Multiplying both sides of (\ref{radau1}) by $x-1$
and using the identity
$$
(x^2-1)L_N'(x) = (N+1)L_{N+1}(x) - (N+1)xL_N(x),
$$
we have
\begin{eqnarray*}
&& N(x-1)(L_N+L_{N-1})(x) \\
&=& (N+1)L_{N+1}(x) - (N+1)xL_N(x) - NL_N(x) + NxL_{N-1}(x).
\end{eqnarray*}
Canceling and collecting the same terms on two sides, we obtain
$$
(2N+1)xL_N(x) - NL_{N-1}(x) = (N+1)L_{N+1}(x),
$$
which is the three-term recurrence relation.

\vskip.1in

{\it Remark 3.1}. We see that zeros of the derivatives of the right (left) Legendre-Radau
polynomials are zeros of the left (right) Legendre-Radau polynomials.

\vskip.1in

\noindent{\bf Lemma 3.2}.
\begin{equation}\label{cheb-radau1}
(T_N+T_{N-1})(x) = (x+1)(U_{N-1}-U_{N-2})(x).
\end{equation}
\begin{equation}\label{cheb-radau2}
(T_N-T_{N-1})(x) = (x-1)(U_{N-1}+U_{N-2})(x).
\end{equation}

Proof: By the definition,
$$
(T_N+T_{N-1})(x) = \cos N\theta + \cos(N-1)\theta = 2\cos(N-\frac{1}{2})\theta \cos\frac{\theta}{2},
$$
$$
(U_{N-1}-U_{N-2})(x) = \frac{\sin N\theta}{\sin\theta} - \frac{\sin(N-1)\theta}{\sin\theta}
= \frac{\cos(N-\frac{1}{2})\theta}{\cos\frac{\theta}{2}}.
$$
Therefore,
\begin{eqnarray*}
(T_N+T_{N-1})(x) &=& 2\cos^2\frac{\theta}{2} (U_{N-1}-U_{N-2})(x) \\
&=& (1+\cos\theta) (U_{N-1}-U_{N-2})(x) = (1+x) (U_{N-1}-U_{N-2})(x).
\end{eqnarray*}
This established (\ref{cheb-radau1}). The proof of (\ref{cheb-radau2}) is similar and hence is omitted.

\vskip.1in

{\it Remark 3.2}. There is an interesting similarity between (\ref{cheb-radau1}) and (\ref{radau1}),
as well as (\ref{cheb-radau2}) and (\ref{radau2}). Note that $T_N'(x)= NU_{N-1}(x)$, and therefore,
zeros of the derivatives of the right (left) Chebyshev-Radau
polynomials are 'almost' zeros of the left (right) Chebyshev-Radau polynomials.

\vskip.1in

\noindent {\bf Theorem 3.1}. Let $u \in {\cal P}_{N+1}$ and $u_N \in {\cal P}_N$
satisfy (\ref{ode}) with the collocation points $x_k$s being the roots of  $L_N$.
Then we have (up to a constant),
\begin{equation}\label{th1}
(u-u_N)(x) = \int_{-1}^x L_N(t) dt = \frac{1}{2N+1} (L_{N+1}-L_{N-1})(x) = \frac{x^2-1}{N(N+1)}L_N'(x).
\end{equation}
Proof: We see that $(u-u_N)' \in {\cal P}_N$. By the definition of the interpolation points,
we have, up to a constant, $(u-u_N)'(x) = L_N(x)$.
Using the initial condition $u_N(-1) = u(-1)$, the conclusion follows by integration.


\vskip.1in

{\it Remark 3.3}. Theorem 3.1 says that when interpolating derivative at
the $N$ Gauss points, the function value approximation is superconvergent at the
$N-1$ interior Lobatto points, i.e., roots of $L_N'$.

\vskip.1in

\noindent {\bf Theorem 3.2}. Let $u \in {\cal P}_{N+1}$ and $u_N \in {\cal P}_N$ satisfy (\ref{ode})
with the collocation points $x_k$s being the roots of $L_N-L_{N-2}$.
Then we have (up to a constant),
\begin{equation}\label{th2}
u(x) - u_N(x) = \frac{(x^2-1)(2N-1)L_{N-1}(x)}{N(N+1)} - \frac{4N-2}{N(N+1)(2N-3)} (L_{N-1}-L_{N-3})(x).
\end{equation}
Proof: From $(u-u_N)' \in {\cal P}_N$, we have (up to a constant), $(u-u_N)'(x) = L_N(x)-L_{N-2}(x)$.
By the identity
$$
[(1-x^2)L_n'(x)]' + n(n+1) L_n(x) = 0,
$$
we have
$$
u'(x) - u_N'(x) = \frac{ [(1-x^2)L_{N-2}'(x)]' }{(N-2)(N-1)} - \frac{ [(1-x^2)L_N'(x)]' }{N(N+1)}.
$$
Integrating with the initial condition and using the identity
$$
(2n+1)L_n(x) = (L_{n+1}-L_{n-1})'(x),
$$
we derive
\begin{eqnarray*}
&& u(x) - u_N(x) = \frac{ (1-x^2)L_{N-2}'(x) }{(N-2)(N-1)} - \frac{ (1-x^2)L_N'(x) }{N(N+1)} \\
&=& \frac{(x^2-1)(2N-1)L_{N-1}(x)}{N(N+1)} + \frac{4N-2}{(N+1)N(N-1)(N-2))} (1-x^2)L_{N-2}'(x),
\end{eqnarray*}
which is the right-hand side of (\ref{th2}) by the identity
$$
\frac{1}{2n + 1} (L_{n+1}(x) - L_{n-1}(x)) = \frac{1}{n(n + 1)}(x^2-1)L_n'(x).
$$

{\it Remark 3.4}. We see that the magnitude of the first term on the right-hand side of (\ref{th2})
is larger than that of the second term by a factor about $N$. Therefore, the function value approximation
reaches its best at roots of $L_{N-1}$. In other words,
Theorem 3.2 says that when interpolating derivative at
the $N$ Lobatto points, function value approximation is superconvergent at the $N-1$ Gauss points.

\vskip.1in

\noindent {\bf Theorem 3.3}. Let $u \in {\cal P}_{N+1}$ and $u_N \in {\cal P}_N$ satisfy (\ref{ode}),
and let the collocation points $x_k$s be the roots of $L_N \mp L_{N-1}$.
Then we have (up to a constant),
\begin{equation}\label{th3}
(u-u_N)(x) = \frac{N^2}{N^2-1} (L_N \pm L_{N-1})(x)(x \mp 1) - \frac{N}{N^2-1}(L_N \mp L_{N-1})(x)(x \pm 1).
\end{equation}

Proof: By (\ref{radau2}), we have (up to a constant),
$$
(u-u_N)'(x) = N(L_N \mp L_{N-1})(x) = (L_N \pm L_{N-1})'(x)(x \mp 1).
$$
Integrating both sides with the initial condition and using (\ref{radau1})-(\ref{radau2}), we have
\begin{eqnarray*}
&& (u-u_N)(x) = \int_{-1}^x (L_N \pm L_{N-1})'(t)(t \mp 1) dt \\
&=& (L_N \pm L_{N-1})(x)(x \mp 1) - \frac{1}{N}(L_N \mp L_{N-1})(x)(x \pm 1)
+ \frac{1}{N} \int_{-1}^x (L_N \mp L_{N-1})(t)dt.
\end{eqnarray*}
Multiplying both sides by $N$ and moving the last term on the right-hand side to the left-hand side,
we obtain (\ref{th3}).

\vskip.1in

{\it Remark 3.5}. We see that the magnitude of the first term on the right-hand side of (\ref{th3})
is larger than that of the second term by a factor $N$.
Therefore, function value approximation reaches its best at the roots of $L_N \pm L_{N-1}$ for (\ref{th3}).
In other words, Theorem 3.3 says that when interpolating derivative at
the $N$ left (right) Radau points, the function value approximation is superconvergent at
the $N$ right (left) Radau points.

\vskip.1in

Now we turn to the Chebyshev polynomials. We need the following two identities:
\begin{equation}\label{cheb-id1}
( \sqrt{1-x^2}T_n'(x) )' + \frac{n^2}{\sqrt{1-x^2}} T_n(x) = 0,
\end{equation}
\begin{equation}\label{cheb-id2}
\frac{1}{2}( T_{n-1}(x) - T_{n+1}(x) ) = (x^2-1) U_{n-1}(x),
\end{equation}

\noindent {\bf Theorem 3.4}. Let $u \in {\cal P}_{N+1}$ and $u_N \in {\cal P}_N$ satisfy (\ref{ode}),
and let the collocation points $x_k$s be the $N$ roots of $T_N$.
Then we have (up to a constant),
\begin{equation}\label{th5}
(u-u_N)(x) = \frac{N}{N^2-1}(x^2-1)U_{N-1}(x) - \frac{xT_N(x)+(-1)^N}{N^2-1}.
\end{equation}

Proof: We have (up to a constant) $(u-u_N)'(x) = T_N(x)$. By (\ref{cheb-id1}), we write
$$
(u-u_N)'(x) = - \frac{\sqrt{1-x^2}}{N^2} ( \sqrt{1-x^2}T_N'(x) )'.
$$
Integrating and using the initial condition, we derive
\begin{eqnarray*}
(u-u_N)(x) &=& \int_{-1}^x T_N(t) dt = - \int_{-1}^x \frac{\sqrt{1-t^2}}{N^2} ( \sqrt{1-t^2}T_N'(t) )' dt \\
&=& - \frac{1-x^2}{N^2} T_N'(x) - \frac{1}{N^2} \int_{-1}^x tT_N'(t) dt \\
&=& \frac{x^2-1}{N} U_{N-1}(x) - \frac{xT_N(x)+(-1)^N}{N^2} + \frac{1}{N^2} \int_{-1}^x T_N(t) dt.
\end{eqnarray*}
Moving the last term to the left and multiplying the resultant by
$\displaystyle \frac{N^2}{N^2-1}$, we obtain (\ref{th5}).

\vskip.1in

{\it Remark 3.6}. We see that the magnitude of the first term on the right-hand side of (\ref{th5})
is larger than that of the second term by a factor $N$. Therefore, the function value approximation
reaches its best at roots of $U_{N-1}$. In other words,
Theorem 3.4 says that when interpolating derivative at
the $N$ Chebyshev points of the first kind, function value approximation is superconvergent at the $N-1$
Chybyshev points of the second kind.

\vskip.1in

\noindent {\bf Theorem 3.5}. Let $u \in {\cal P}_{N+1}$ and $u_N \in {\cal P}_N$ satisfy (\ref{ode}),
and let the collocation points $x_k$s be $\pm 1$ plus the $N-2$ roots of $U_{N-2}$.
Then we have (up to a constant),
\begin{eqnarray}\label{th6}
(u-u_N)(x)
&=& \frac{N(1-x^2)}{N^2-1} T_{N-1}(x) - \frac{N(1-x^2)}{(N^2-1)(N-2)} U_{N-3}(x)
+ \frac{xT_N(x) + (-1)^N}{2(N^2-1)}   \nonumber \\
&& + \frac{N^2(xT_{N-2}(x) + (-1)^N)}{2(N-2)^2(N^2-1)} + \frac{N-1}{(N^2-1)(N-2)^2} \int_{-1}^x T_{N-2}(t)dt.
\end{eqnarray}

Proof: We have (up to a constant), $(u-u_N)'(x) = (x^2-1)U_{N-2}(x)$.
By (\ref{cheb-id2}) and (\ref{cheb-id1}), we write
\begin{eqnarray*}
&& (u-u_N)'(x) = \frac{1}{2} ( T_{N-2}(x) - T_N(x) ) \\
&=& \frac{\sqrt{1-x^2}}{2} \left( \frac{1}{N^2}( \sqrt{1-x^2}T_N'(x) )'
- \frac{1}{(N-2)^2}( \sqrt{1-x^2}T_{N-2}'(x) )' \right).
\end{eqnarray*}
Integrating and using the initial condition, we derive
\begin{eqnarray*}
&& (u-u_N)(x) = \frac{1}{2} \int_{-1}^x ( T_{N-2}(t) - T_N(t) ) dt \\
&=& \frac{1-x^2}{2N^2} T_N'(x) + \frac{1}{2N^2} \int_{-1}^x tT_N'(t) dt
- \frac{1-x^2}{2(N-2)^2}T_{N-2}'(x) - \frac{1}{2(N-2)^2} \int_{-1}^x tT_{N-2}'(t) dt \\
&=& \frac{1-x^2}{2N} (U_{N-1}(x) - U_{N-3}(x)) - \frac{1-x^2}{N(N-2)} U_{N-3}(x) + \frac{xT_N(x) + (-1)^N}{2N^2} \\
&& + \frac{xT_{N-2}(x) + (-1)^N}{2(N-2)^2} + \frac{N-1}{N^2(N-2)^2} \int_{-1}^x T_{N-2}(t)dt
+ \frac{1}{2N^2} \int_{-1}^x (T_{N-2}(t)-T_N(t)) dt.
\end{eqnarray*}
Moving the last term to the left-hand side and multiplying the resultant
by $\displaystyle \frac{N^2}{N^2-1}$, we derive (\ref{th6}) when replacing
$U_{N-1}-U_{N-3}$ by $2T_{N-1}$.

\vskip.1in

{\it Remark 3.7}. We see that the magnitude of the first term on the right-hand side of (\ref{th6})
is larger than that of the rest terms by a factor about $N$. Therefore, $(u-u_N)(x)$
reaches its best at roots of $T_{N-1}$. In other words,
Theorem 3.5 says that when interpolating derivative at $\pm 1$ plus
the $N-2$ Chebyshev points of the second kind, function value approximation is superconvergent at the $N-1$
Chebyshev points of the first kind.

\vskip.1in

\noindent {\bf Theorem 3.6}. Let $u \in {\cal P}_{N+1}$ and $u_N \in {\cal P}_N$ satisfy (\ref{ode}),
and let the collocation points $x_k$s be the roots of $T_N \pm T_{N-1}$.
Then we have (up to a constant),
\begin{eqnarray}\label{th7}
(u-u_N)(x)
&=& \frac{N(x \pm 1)}{N^2-1} (T_N \mp T_{N-1})(x) \pm \frac{N(x^2-1)}{(N^2-1)(N-1)} U_{N-2}(x)
\mp \frac{xT_N(x) + (-1)^N}{N^2-1} \nonumber\\
&& \mp \frac{N^2(xT_{N-1}(x) + (-1)^{N-1})}{(N^2-1)(N-1)^2}
\pm \frac{2N-1}{(N^2-1)(N-1)^2} \int_{-1}^x T_{N-1}(t) dt.
\end{eqnarray}

Proof: The proof is similar to that of Theorem 3.5 and hence is omitted.


\vskip.1in

{\it Remark 3.8}. We see that the magnitude of the first term
on the right-hand side of (\ref{th7})
is larger than that of the rest terms by a factor about $N$.
Therefore, the function value approximation
reaches its best at the roots of $T_N \mp T_{N-1}$ for (\ref{th7}).
In other words, Theorems 3.6 says that when interpolating derivative at
the $N$ left (right) Chebyshev-Radau points, function value approximation is superconvergent at
the $N$ right (left) Chebyshev-Radau points.

\bigskip


\setcounter{equation}{0}
\setcounter{section}{4}

\noindent{\bf 4. Numerical Tests}

In this section, we perform numerical tests on two typical analytic functions.

\begin{figure} 
        \centering
        \includegraphics[scale=0.5]{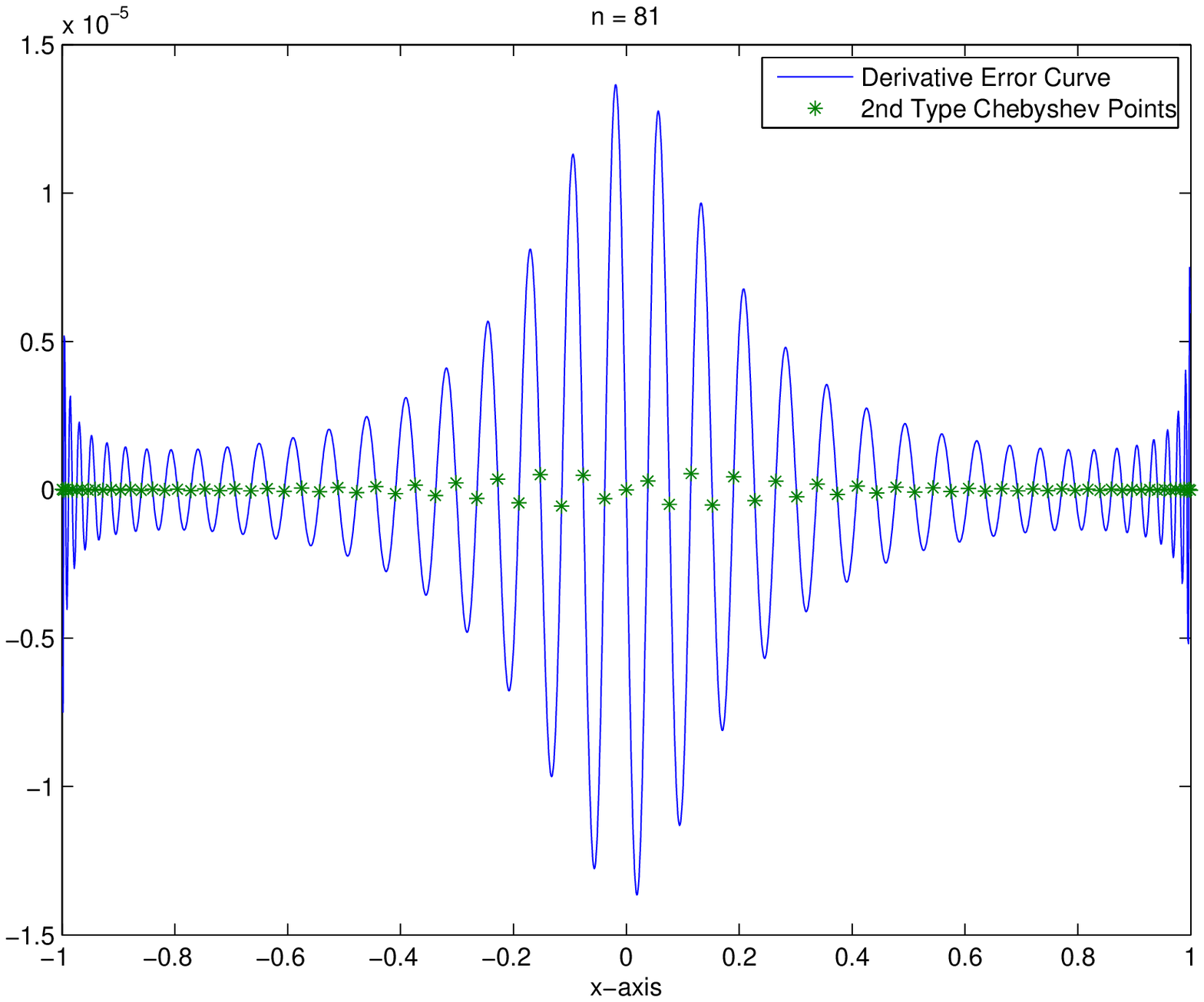}
    \caption{Interpolation at the Chebyshev Points - Example 1}
    \label{ch1}
\end{figure}

\begin{figure} 
        \centering
        \includegraphics[scale=0.5]{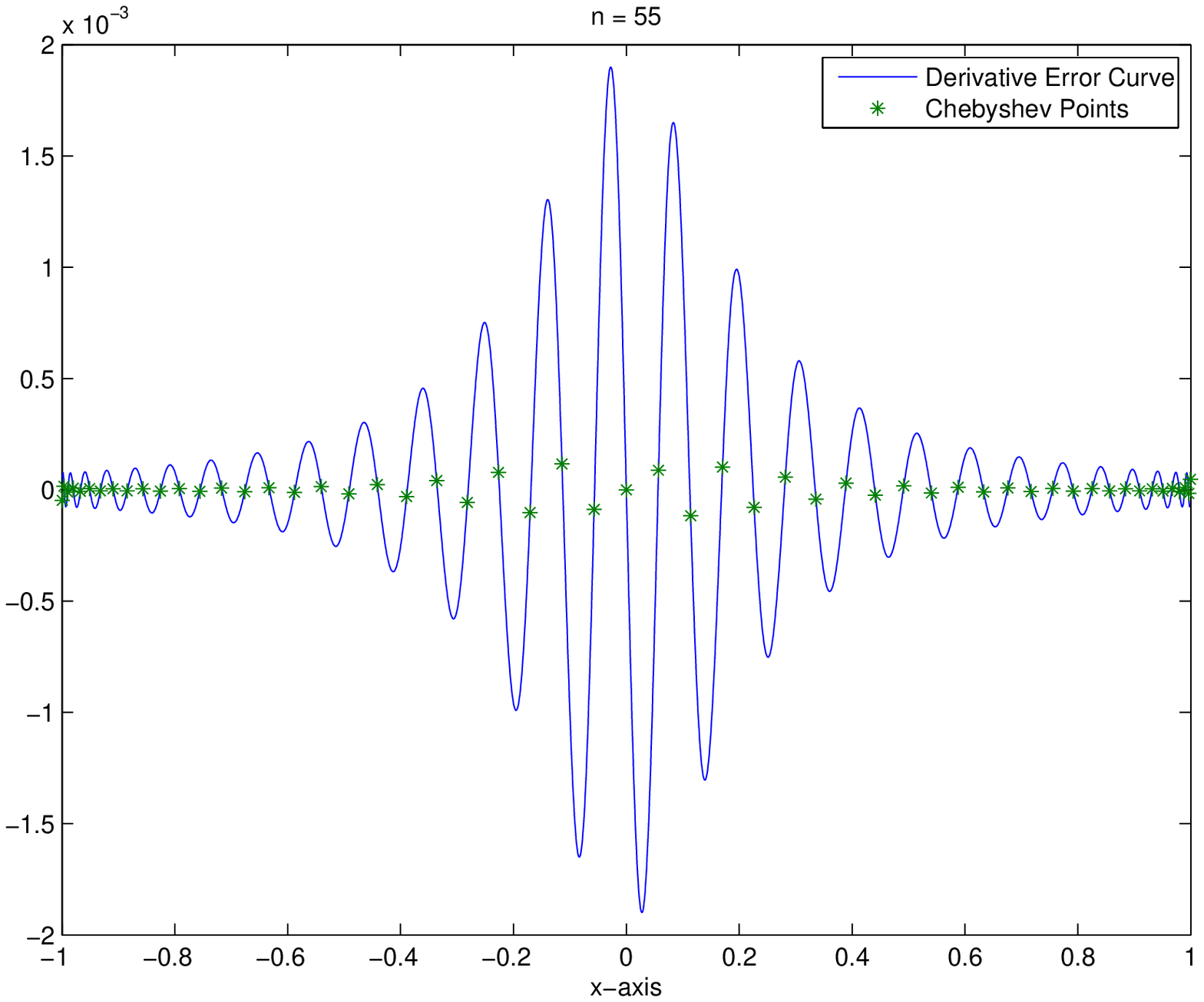}
    \caption{Interpolation at the 2nd Type Chebyshev Points - Example 1}
    \label{ch-lobatto1}
\end{figure}

\begin{figure} 
        \centering
        \includegraphics[scale=0.5]{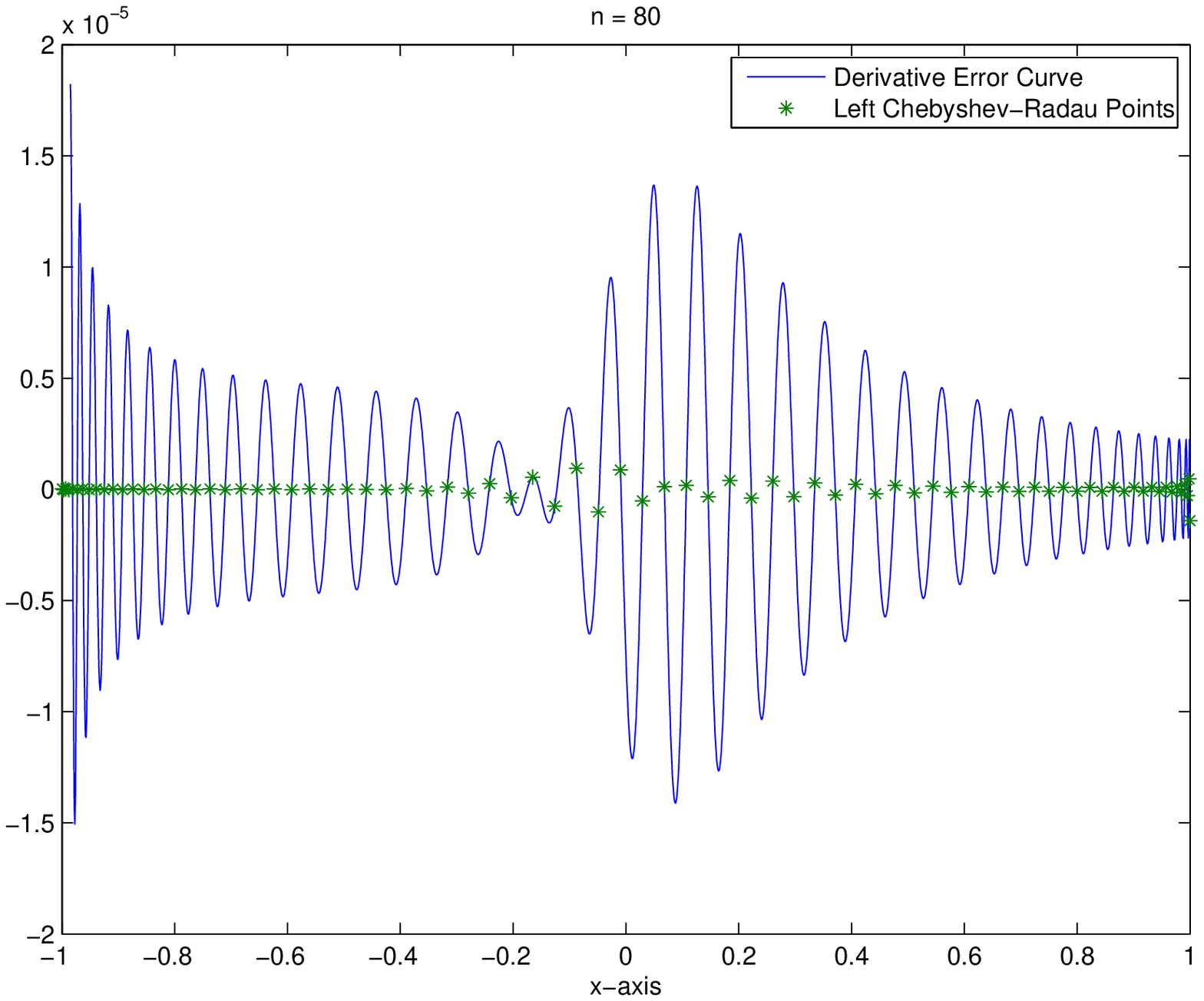}
    \caption{Interpolation at the Right Chebyshev-Radau Points - Example 1}
    \label{ch-radau1}
\end{figure}

{\it Example 1}. $f(x) = (1+25x^2)^{-1}$. This is the well known  Runge's example \cite{runge}.
Derivative errors of its interpolants at the Chebyshev points, the 2nd type Chebyshev points
(or Chebyshev-Lobatto points),
and the right-Chebyshev-Radau points are depicted in Figure \ref{ch1}, Figure \ref{ch-lobatto1}, and
Figure \ref{ch-radau1}, respectively. Derivative superconvergent points are marked by '*'.
We see that the errors at the superconvergent points are significantly smaller (by a magnitude)
than the maximum error just as Theorems 2.1 - 2.3 predicted.

\begin{figure} 
        \centering
        \includegraphics[scale=0.5]{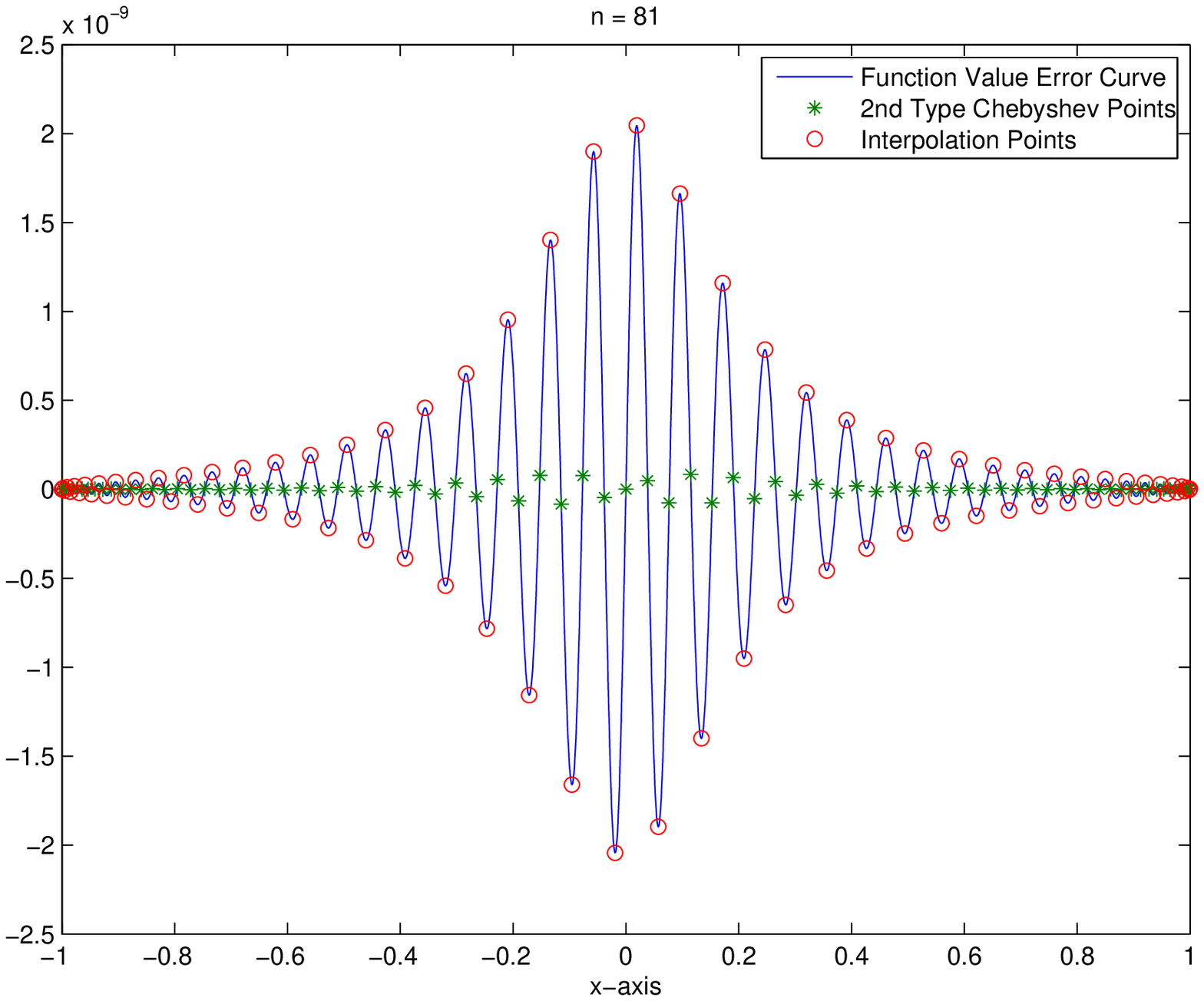}
    \caption{Derivative Interpolation at the Chebyshev Points - Example 1}
    \label{ch-diff1}
\end{figure}

Next, we solve the following initial value problem
$$
u'(x) = \frac{1}{1+25x^2}, \quad u(-1) = \frac{1}{26}
$$
by collocating at the Chebyshev points as in (4.1). From Figure \ref{ch-diff1},
we see that when interpolating the derivative at the Chebyshev points,
the function value of the interpolant converges much faster at the interior 2nd type Chebyshev points
marked by '*', which is in consistent with Theorem 4.1. We also see that the function value
approximation has maximum error at the derivative interpolation points marked by '$\circ$'.

Analytic function $f(z) = (1+25x^2)^{-1}$ has two single poles at $\pm i/5$ and it is
straightforward to calculate $\rho = (\sqrt{5^2+1}+1)/5 \approx 1.2198$. Therefore, we expect a slow
convergence as we have observed form Figure \ref{ch1} - Figure \ref{ch-diff1}.

\vskip.1in

\begin{figure} 
        \centering
        \includegraphics[scale=0.5]{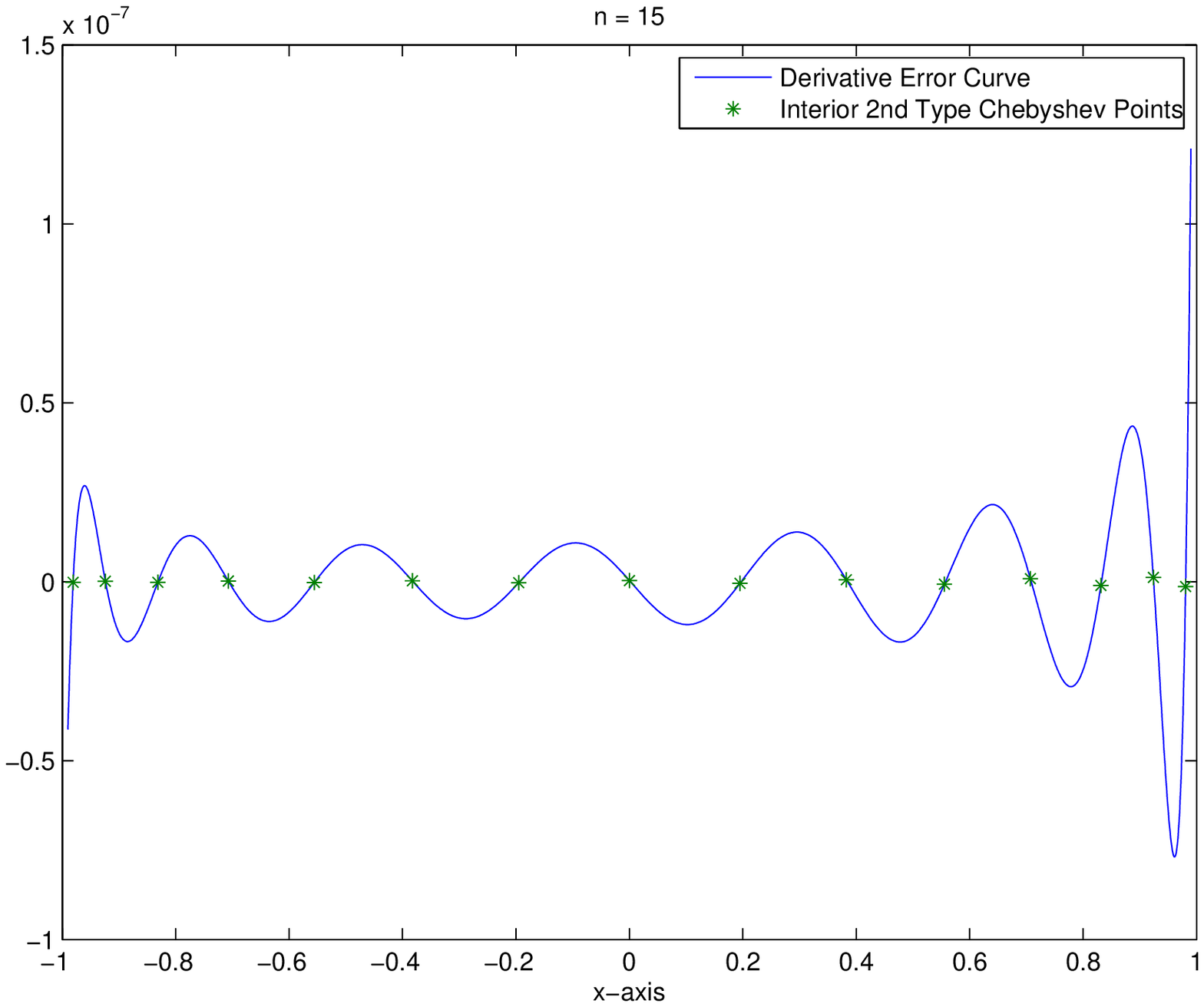}
    \caption{Interpolation at the Chebyshev Points - Example 2}
    \label{ch2}
\end{figure}

\begin{figure} 
        \centering
        \includegraphics[scale=0.5]{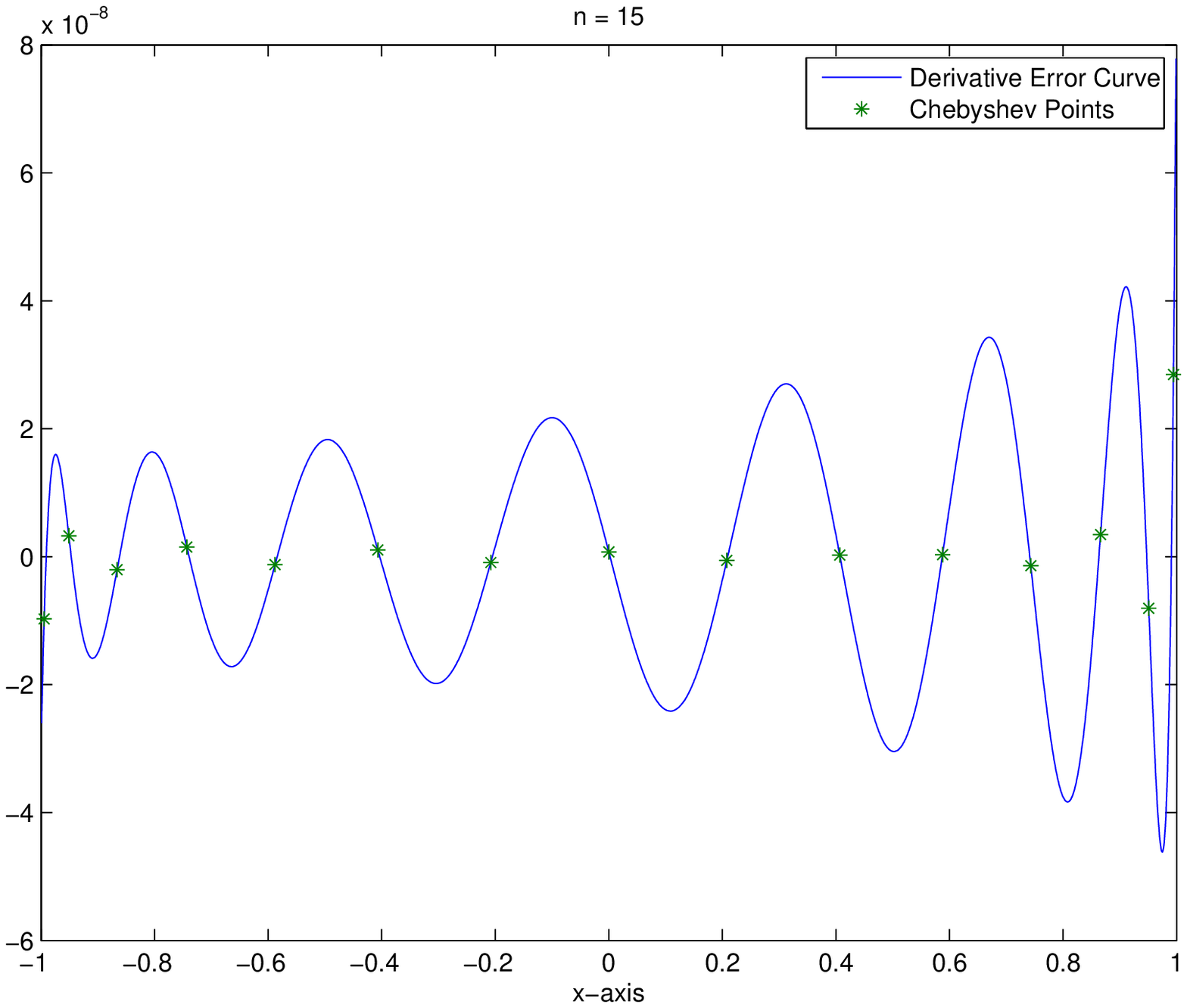}
    \caption{Interpolation at the 2nd Type Chebyshev Points - Example 2}
    \label{ch-lobatto2}
\end{figure}

\begin{figure} 
        \centering
        \includegraphics[scale=0.5]{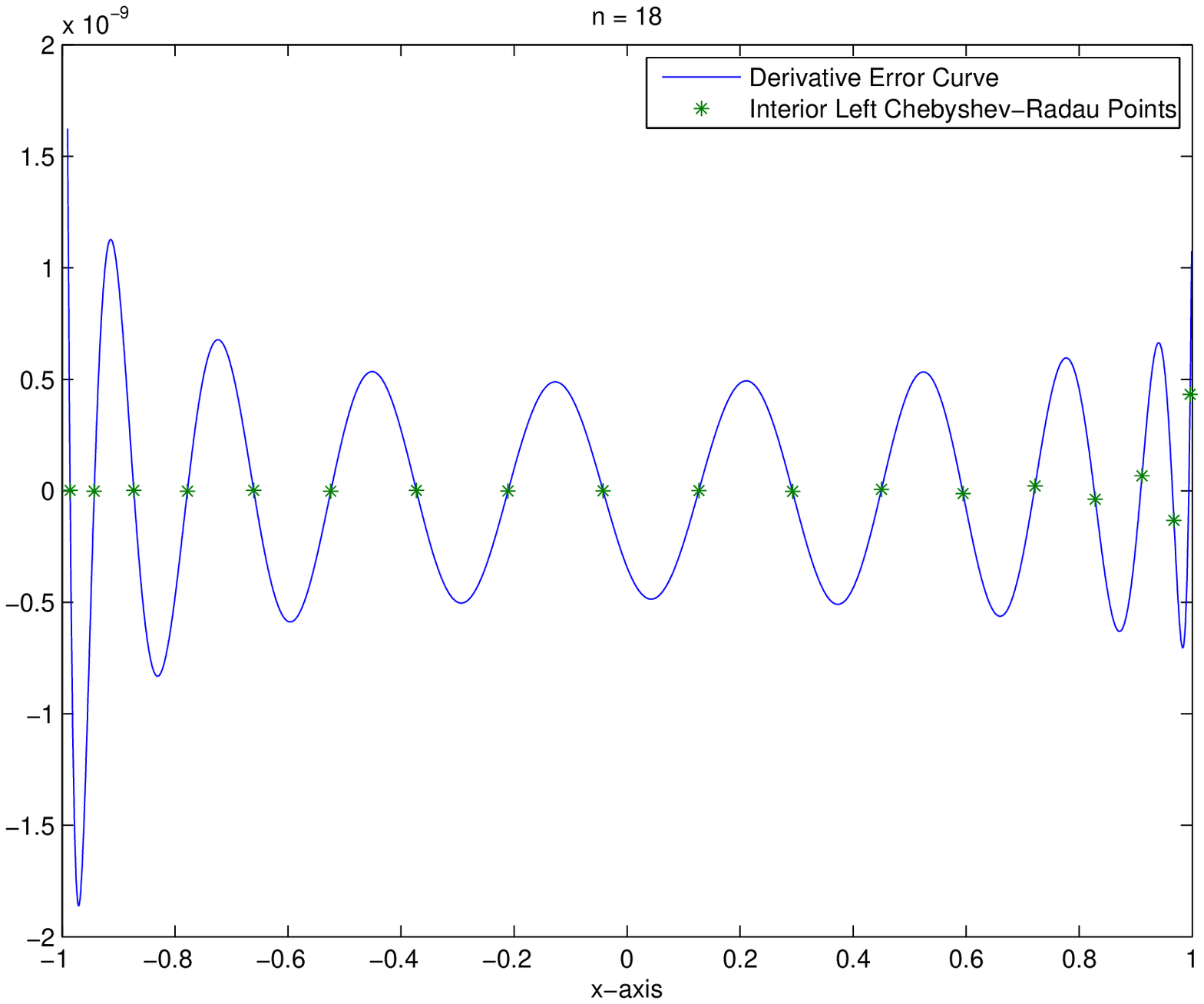}
    \caption{Interpolation at the Right Chebyshev-Radau Points - Example 2}
    \label{ch-radau2}
\end{figure}

\begin{figure} 
        \centering
        \includegraphics[scale=0.5]{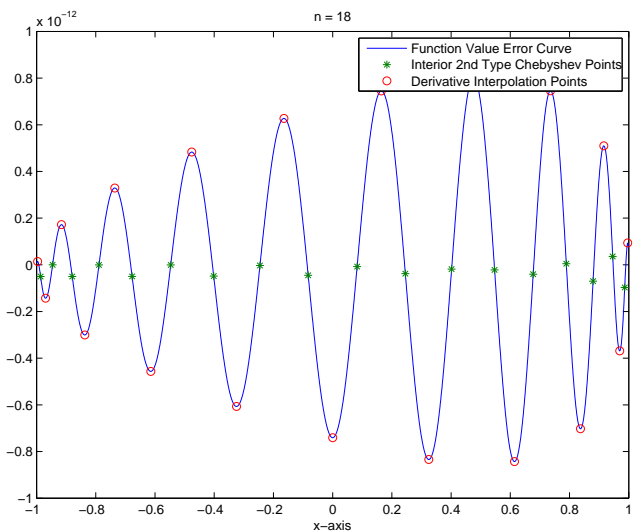}
    \caption{Derivative Interpolation at the Chebyshev Points - Example 2}
    \label{ch-diff2}
\end{figure}

{\it Example 2}. $f(x) = (2-x)^{-1}$, an analytic function which has a simple pole at $x=2$
and $\rho = 2 + \sqrt{2^2-1} \approx 3.7321$. Therefore, we expect much faster convergence
compared with {\it Example 1}. It is indeed the case. We plot counterparts graphs
in Figure \ref{ch2} -- Figure \ref{ch-diff2}, and we see that high accuracy is achieved with
relative very low polynomial degree $n$. We observe similar superconvergence phenomena as in
{\it Example 1}.

\vskip.2in

{\bf Conclusion Remarks}. The results in Section 2 can be extended to all
Legendre based polynomials as in Section 3.
The reason for using the Chebyshev based polynomials is that they have simple
(trigonometric functions) expressions.

Extension of the result to more general Jacobi type polynomials is feasible.
However, the analysis would be much more involved.

\medskip
\begin{thebibliography} {99}


\bibitem{ainsworth-oden} M. Ainsworth and J.T. Oden,
{\it A Posteriori Error Estimation in Finite Element Analysis},
Wiley Interscience, New York, 2000.

\bibitem{Babenko} K.I. Babenko,
On the best approximation of a class of analytic functions,
Izv. 22 (1958), 631-640.

\bibitem{babuska-strouboulis}
I. Babu\v{s}ka and T. Strouboulis,
{\it The Finite Element Method and its Reliability},
Oxford University Press, London, 2001.

\bibitem{babuska-sug}
I. Babu\v{s}ka, T. Strouboulis, C.S. Upadhyay, and S.K. Gangaraj,
Computer-based proof of the existence of superconvergence points in the finite element method:
superconvergence of the derivatives in finite element solutions of Laplace's, Poisson's, and the elasticity equations, {\it Numer. Meth. PDEs.} {\bf 12} (1996), 347-392.

\bibitem{bernardi-maday} C. Bernardi and Y. Maday,
Spectral Methods, in: P.G. Ciarlet and J.-L. Lions (eds.),
{\it Handbook of Numerical Analysis} Vol. V, North-Holland (1997), 209-485.

\bibitem{bernstein} S.N. Bernstein,
Sur l'ordre de la meilleure approximation des foncions
continues par des polynomes de degr\'{e} donn\'{e},
M\'{e}moires publi\'{e}s par la class des sci. Acad. de Belgique (2) 4 (1912), 1-103.



\bibitem{boyd} J.P. Boyd,
{\it Chebyshev and Fourier Spectral Methods}, 2nd edition, Dover, New York, 2001.

\bibitem{canuto-hqz} C. Canuto, M.Y. Hussaini, A. Quarteroni, and T.A. Zang,
{\it Spectral Methods in Fluid Dynamics}, Springer-Verlag, New York, 1988.

\bibitem{canuto-hqz1} C. Canuto, M.Y. Hussaini, A. Quarteroni, and T.A. Zang,
{\it Spectral Methods: Fundamentals in Single Domains},
Springer-Verlag, New York, 2006.

\bibitem{canuto-hqz2} C. Canuto, M.Y. Hussaini, A. Quarteroni, and T.A. Zang,
{\it Spectral Methods: Evolution to Complex Geometries and Applications to Fluid Dynamics},
Springer-Verlag, New York, 2007.


\bibitem{chen} C.M. Chen,
{\it Structure Theory of Superconvergence of Finite Elements} (in Chinese),
Hunan Science and Technology Press, China, 2001.

\bibitem{chen-huang} C.M. Chen and Y.Q. Huang,
{\it High Accuracy Theory of Finite Element Methods} (in Chinese),
Hunan Science and Technology Press, China, 1995.

\bibitem{davis} P.J. Davis, {\it Interpolation and Approximation},
Dover Publications Inc., New York, 1975.



\bibitem{davis-rabinowitz} P.J. Davis and P. Rabinowitz,
{\it Methods of Numerical Integration}, 2nd edition,
Academic Press, Boston, 1984.

\bibitem{day-romero}  David Day and Louis Romero, Roots of polynomials expressed in terms of orthogonal polynomials, {\it SIAM J. Numer. Anal.} 43-5 (2006), 1969-1987.

\bibitem{deville-fm} M.O. Deville, P.F. Fischer, and E.H. Mund,
{\it High-Order Methods for Incompressible Fluid Flow},
Cambridge University Press, 2002.

\bibitem{fornberg} B. Fornberg.
{\it A Practical Guide to Pseudospectral Methods},
Cambridge University Press, 1996.

\bibitem{funaro} D. Funaro.
{\it Polynomial Approximations of Differential Equations},
Springer-Verlag, New York, 1992.

\bibitem{gottlieb-orszag} D. Gottlieb and T.A. Orszag,
{\it Numerical Analysis of Spectral Methods: Theory and Applications},
SIAM, Philadelphia, 1977.


\bibitem{gottlieb-shu} D. Gottlieb and C.-W. Shu,
On the Gibbs phenomenon and its resolution,
{\it SIAM Rev.} 39-4 (1997), 644-668.

\bibitem{guo-babuska1} B.Q. Guo and I. Babu\v{s}ka,
The h-p version of the finite element method, Part 1:
the basic approximation results, {\it Comp. Mech.} 1 (1986), 22-41.

\bibitem{guo-babuska2} B.Q. Guo and I. Babu\v{s}ka,
The h-p version of the finite element method, Part 2:
the general results and application, {\it Comp. Mech.} 1 (1986), 203-220.


\bibitem{guo} B.Y. Guo. {\it Spectral Methods and Their Applications},
World Scientific Publishing Company, Beijing, 1998.

\bibitem{hesthaven-gg} J.S. Hesthaven, S. Gottlieb, and D. Gottlieb,
{\it Spectral Methods for Time-Dependent Problems},
Cambridge University Press, 2007.

\bibitem{karniadakis-sherwin} G.E. Karniadakis and S.J. Sherwin,
{\it Spectral/hp Element Methods for CFD},
Oxford University Press, New York, 1998.

\bibitem{krizek-n}
M. K\v{r}\'{\i}\v{z}ek, and P. Neittaanm\"{a}ki,
On superconvergence techniques,
{\it Acta Appl. Math.} {\bf 9} (1987), 175-198.

\bibitem{krizek-ns}
M. K\v{r}\'{\i}\v{z}ek, P. Neittaanm\"{a}ki, R. and Stenberg (Eds.),
{\it Finite Element Methods: Superconvergence, Post-processing, and A Posteriori Estimates},
Lecture Notes in Pure and Applied Mathematics Series Vol.\ 196, Marcel Dekker Inc., New York, 1997.




\bibitem{lin-lin} Q. Lin and J. Lin,
{\it Fininte Element Methods: Accuracy and Improvement},
Mathematics Monograph Series Vol.\ 1, Science Press, Beijing, China, 2006.

\bibitem{lin-yan} Q. Lin and N. Yan,
{\it Construction and Analysis of High Efficient Finite Elements} (in Chinese),
Hebei University Press, China, 1996.

\bibitem{lorentz} G.G. Lorentz, {\it Approximation of Functions},
AMS Chelsea Publishing, 1966.

\bibitem{mason-handscomb} J.C. Mason and D.C. Handscomb,
{\it Chebyshev Polynomials},
Chapman \& Hall/CRC, Boca Raton, 2003.



\bibitem{peyret} R. Peyret,
{\it Specrtral Methods for Incompressible Viscous Flow},
Springer, New York, 2002.



\bibitem{phillips} G.M. Phillips,
{\it Interpolation and Approximation by Polynomials},
Springer, New York, 2003.

\bibitem{pozrikidis} C. Pozrikidis,
{\it Introduction to Finite and Spectral Element Methods Using {\tt Matlab}},
Chapman \& Hall/CRC, 2005.

\bibitem{runge} Carl Runge,
\"{U}ber empirische Funktionen und die Interpolation zwischen \"{a}quidistanten Ordinaten,
{\it Zeitschrift f\"{u}r Mathematik und Physik} 46 (1901), 224-243.

\bibitem{sansone} G. Sansone, {\it Orthogonal Functions},
Dover, New York, 1991.



\bibitem{reddy-weideman} S.C. Reddy and J.A.C. Weideman,
The accuracy of the Chebyshev differencing method for analytic functions,
{\it SIAM J. Numer. Anal.} 42-2 (2005), 2176-2187.

\bibitem{rivlin} T.J. Rivlin,
{\it An Introduction to the Approximation of Functions}, Dover, New York, 1969.

\bibitem{schatz-sloan-wahlbin}
A.H. Schatz, I.H. Sloan, and L.B. Wahlbin,
Superconvergence in finite element methods and meshes which are symmetric with respect to a point,
{\it SIAM J.\ Numer. Anal.} {\bf 33} (1996), 505-521.

\bibitem{schwab} C. Schwab, {\it p- and hp- Finite Element Methods},
Oxford University Press, 1998.


\bibitem{shen-tang} J. Shen and T. Tang,
{\it Spectral and High-Order Methods with Applications},
Science Press of China, Beijing, 2006.

\bibitem{shen-tang-wang} J. Shen, T. Tang, and L.-L. Wang.
{\it Spectral Methods: Algorithms, Analysis and Applications}, Springer, 2011.


\bibitem{szabo-babuska} B. Szab\'{o} and I. Babu\v{s}ka, {\it Finite Element Analysis},
John Wiley \& Sons, Inc., New York, 1991.

\bibitem{szego} G. Szeg\"{o}, {\it Orthogonal Polynomials},
4th edition, AMS Colloq. Public. Vol.23, 1975.

\bibitem{tadmor} E. Tadmor,
The exponential accuracy of Fourier and Chebyshev differencing methods,
{\it SIAM J. Numer. Anal.} 23-1 (1986), 1-10.

\bibitem{tang-xu} T. Tang and J. Xu (eds.),
{\it Adaptive Computations: Theory and Algorithms},
Mathematics Monograph Series 6, Science Publisher, 2007.

\bibitem{trefethen} L.N. Trefethen, {\it Spectral Methods in {\tt Matlab}},
SIAM, 2000.

\bibitem{wahlbin} L.B. Wahlbin,
{\it Superconvergence in Galerkin Finite Element Methods},
Lecture Notes in Mathematics Vol.\ 1605, Springer, Berlin, 1995.

\bibitem{wang-xie-zhang} L. Wang, Z. Xie, and Z. Zhang,
Super-geometric convergence of spectral element method for eigenvalue problems with jump coefficients,
{\it Journal of Computational Mathematics} 28-3 (2010), 418-428.

\bibitem{zhang1} Z. Zhang,
Superconvergence of Spectral collocation and p-version methods in one dimensional problems,
{\it Math. Comp.} 74 (2005), 1621-1636.

\bibitem{zhang2} Z. Zhang,
Superconvergence of a Chebyshev spectral collocation method,
{\it J. Sci. Comp.} 34 (2008), 237-246.

\bibitem{zhu-lin} Q.D. Zhu and Q. Lin,
{\it Superconvergence Theory of the Finite Element Method} (in Chinese),
Hunan Science Press, China, 1989.

\end{thebibliography}
\end{document}